\documentclass{article}

\usepackage{srcltx}
\usepackage{graphicx,color,amssymb,latexsym,amsmath,amsfonts}
\usepackage{mathrsfs}

\RequirePackage{amsthm,amsmath,amssymb}
\RequirePackage[numbers]{natbib}

\numberwithin{equation}{section}
\theoremstyle{plain}
\newtheorem{theorem}{Theorem}[section]

\newtheorem{lemma}{Lemma}[section]

\title{Chi-squared test for  hypothesis testing  of homogeneity}
\author{Mikhail Ermakov}

\begin{document}
\global\long\def\zb{\boldsymbol{z}}
\global\long\def\ub{\boldsymbol{u}}
\global\long\def\vb{\boldsymbol{v}}
\global\long\def\wb{\boldsymbol{w}}
\global\long\def\tb{\boldsymbol{t}}
\global\long\def\eb{\boldsymbol{e}}
\global\long\def\sb{\boldsymbol{s}}

\maketitle

\footnote{This Research has been supported RFFI Grant   20-01-00273.}

key words: goodness of fit tests, consistency, chi-squared test, maxisets.

\begin{abstract} We provide necessary and sufficient conditions of uniform consistency of nonparametric sets of alternatives of  chi-squared test for testing of hypothesis of homogeneity. The number of cells of chi-squared test increases with sample size growth. Nonparametric sets of alternatives can be defined both in terms of densities and distribution functions.
\end{abstract}

\section{Introduction}
For goodness-of-fit testing chi-squared tests with increasing number of cells with growth of sample size are comprehensively studied \cite{bar,chib,er97,er20,iv,ing87,man,mor,vaart}.

Let $X_1,\ldots,X_n$ be sample of i.i.d.r.v.'s with values on interval $[0,1]$ and having c.d.f.  $F_n$.  Denote $\hat F_n$ -- empirical c.d.f. of sample.
Denote $\Im $--set of all distribution functions. Denote $F_0$ --c.d.f. of uniform distribution on interval  $[0,1]$. Goodness-of-fit testing we  consider as problem of testing hypothesis $\mathbb{H}_0\,:\,F_n=F_0$ versus alternatives $\mathbb{H}_n\,:\,F_n \in \Psi_n \subset \Im$, where $\Psi_n$ is some nonparametric set of alternatives.  

Denote $T_n(\hat F_n)$ -- test statistics of chi-squared tests and  $T_n(F)$, $F \in \Im$, -- functionals generating test statistics $T_n(\hat F_n)$

For goodness-of-fit testing we show in \cite{er97} that sequences of chi-squared tests having increasing number of cells with growth of sample size are uniformly consistent on sets of alternatives  $\Im(b_n) =\{ F\,:\, T_n(F) > b_n, F \in \Im\,\}$, where  $b_n >0$ is sequence of constant depending on number of cells and sample size $n$. Thus sequence of sets of alternatives $\Omega_n \subset \Im$ is uniformly consistent, if and only if,  $\Omega_n \subset \Im(b_n)$ with sequence of numbers $b_n$ satisfying certain conditions. In \cite{er20} we described   all uniformly consistent sequences of alternatives defined in terms of densities if cells of chi-squared test have  equal length   and number of cells growth with increasing sample size.

Paper goal is to  explore uniform consistency of chi-squared tests having increasing number of cells with growth of sample size for testing of hypothesis homogeneity. The goal is to describe all uniformly consistent sequences of alternatives  defined in terms of distribution functions or densities. The problem is more difficult than for goodness -of-fit testing  \cite{er97,er20}.  For hypothesis testing of homogeneity the answer depends on distribution functions of two samples. Note that problem of hypothesis testing of homogeneity has been  intensively studied in recent papers
 \cite{from12,from13, gr12a,gr12b,tong}.

  Let the interval $ [0,1] $ be divided into $ m = m_n $ subintervals
 \begin{equation*}
 I_{nj}= [e_{nj},e_{n,j+1}), \quad p_{nj} = e_{n,j+1} -e_{nj} > 0, \quad e_{n0} = 0,\quad e_{nm}=1,
  \end{equation*}
 $ 1 \le j \le m= m_n$, where  $m_n \to \infty$ as $n \to \infty$. Functional $T_n$ generating chi-squared test statistics for goodness-of-fit testing equals
 \begin{equation*}
 T_n(F-F_0) = n \sum_{j=1}^m \frac{(r_{nj}  - p_{nj})^2}{p_{nj}},
  \end{equation*}
 where $r_{nj}  = F(e_{nj}) - F(e_{n,j-1})$ for all $1 \le j \le m_n$ and $F_0(x)= x$,  $x \in [0,1]$.

Then $T_n(\hat F_n - F_0)$ is chi-squared test statistics.

For test $K_n$ we denote $\alpha(K_n)$ -- its type I error probability and $\beta(K_n,F_n)$-- its type II error probability for alternative  $F_n$.

 Let $S_n$ be sequence of test statistics. We say that sequence of sets of alternatives  $\Psi_n \subset \Im$ is uniformly consistent for test statistics $S_n$, if for tests $K_n$ generated test statistics  $S_n$ such that $\alpha(K_n) = \alpha(1 + o(1))$, $0 <\alpha <1$, we have
  \begin{equation*}
 \limsup_{n \to \infty} \sup_{F \in \Psi_n} \beta(K_n,F) <\,1 - \alpha.
 \end{equation*}
Similar notation and terminology we shall use for problem of testing of hypothesis if homogeneity. 
  As mentioned, for goodness-of-fit testing chi-squared test is uniformly consistent for sets of alternatives 
$\Im(b_n)$.
Moreover \cite{er97}, for any sequence of simple hypothesis $F_n \in \Im$ for type II error probabilities  $\beta(K_n,F_n)$ of tests $K_n$, $\alpha(K_n) = \alpha(1 + o(1))$, $0 <\alpha <1$,  generated test statistics $T_n(\hat F_n - F_0)$,   we have
\begin{equation}\label{uuu}
\beta(K_n,F_n) = \Phi(x_\alpha - 2^{-1/2}m_n^{-1/2}T_n(F_n-F_0)) +o(1)
\end{equation}
as $n \to \infty$. Here  $x_{\alpha}$  is defined by equation $\alpha = 1 - \Phi(x_\alpha)$, where $\Phi(x)= \frac{1}{\sqrt{2\pi}}\int_{-\infty}^x \exp\{-2\,t^2/2\}\,d\,t$ -- standard normal distribution function, $x \in \mathbb{R}^1$.

Such an asymptotic of type II error probabilities and asymptotic minimaxity of chi-squared tests \cite{er97} substantiates the reasoning for using the method of distances in nonparametric hypothesis testing in relation to chi-squared tests. 

In paper we establish similar results for testing of hypothesis of homogeneity with sets of alternatives generated with differences of distribution functions of  two samples.   We suppose additionally that $\mathbb{L}_2$-norm of densities of one of distribution functions  $F_n$ are bounded some constant. It turns out that uniform consistency of sets of alternatives is given by the  value of functional  $T_n$ defined on differences of distribution functions of these two samples. This allows to extend on this setup the results of   \cite{er20}  on necessary and sufficient conditions of uniform consistency of sets of alternatives defined in terms of densities.

 We use letters $c$ and $C$ as a generic notation for positive constants.   Denote $[a]$ whole part of real number $a$. For any two sequences of positive real numbers $a_n$ and $b_n$,  $a_n \asymp b_n$ implies $c < a_n/b_n < C$ for all $n$ and $a_n = o(b_n)$ implies $a_n/b_n \to 0$ as $n \to \infty$. For any complex number $z$ denote $\bar z$ complex conjugate number.

 \section{Main results}

\subsection{Setup} In comparison with goodness-of-fit-testing the problem more difficult. We have two samples  $X_1,\ldots,X_n$ and $Y_1, \ldots, Y_{l_n}$ of i.i.d.r.v's taking values on interval  $[0,1]$ and having distribution functions  $F_n$  and   $G_{l_n}$ respectively. Thus the criterion of uniform consistency has to be sought in terms of differences $G_n - F_n$ and nuisance parameter $F_n$ or $G_n$.

Denote $\Im\times\Im$ -- set of all pairs of distribution functions  $(F,G)$.

On the set $\Im\times \Im$ we define functional
\begin{equation*}
 T_{1n}(F- G) = n m \sum_{j=1}^m (r_{nj}  - s_{nj})^2,\qquad (F,G) \in \Im \times \Im,
  \end{equation*}
 where $s_{nj}  = G(e_{nj}) - G(e_{n,j-1})$ for all  $1 \le j \le m_n$.

 Denote $\hat G_{l_n}(x)$ -- empirical distribution function of sample   $Y_1, \ldots, Y_{l_n}$.

Denote $a_n = \frac{n}{l_n}$ and suppose that  $0 < c < a_n < C < \infty$.

Chi-squared test statistics has the following form 
\begin{equation*}
 T_{1n}(\hat F_n-\hat G_{l_n}) = n m \sum_{j=1}^m (\hat r_{nj}  - \hat s_{nj})^2,
  \end{equation*}
 where $\hat s_{nj}  = \hat G_{l_n}(e_{nj}) - \hat G_{l_n}(e_{n,j-1})$ for all  $1 \le j \le m_n$.

Note that $\mathbf{E}[ T_{1n}(\hat F_n-\hat G_{l_n})]$ depends only on difference of distribution functions $F_n - G_{l_n}$ and we do not need to add additional estimates of addendums to test statistics \cite{and,er97,from12,from13}, related with dependence on distribution functions $F_n$ and $G_{l_n}$.

 For test statistics
  \begin{equation*}
 T_{2n}(\hat F_n-\hat G_{l_n}) = n \sum_{j=1}^m g_{nj}\frac{(\hat r_{nj}  - \hat s_{nj})^2}{p_{nj}},
  \end{equation*}
 generated functional
 \begin{equation*}
 T_{2n}(F- G) = n \sum_{j=1}^m g_{nj}\frac{(r_{nj}  - s_{nj})^2}{p_{nj}},\quad 0<c < g_{jn}< C<\infty
  \end{equation*}
  a separate theorem is proved.

  Proof is provided for test statistics 
   \begin{equation*}
 T_n(\hat F_n-\hat G_{l_n}) = n \sum_{j=1}^m \frac{(\hat r_{nj}  - \hat s_{nj})^2}{p_{nj}},
  \end{equation*}
generated by functional
 \begin{equation*}
 T_{n}(F- G) = n \sum_{j=1}^m \frac{(r_{nj}  - s_{nj})^2}{p_{nj}}.
   \end{equation*}
  For test statistics $ T_{2n}(\hat F_n-\hat G_{l_n})$ the reasoning are almost the same and therefore the differences are not indicated.

We suppose that nuisance parameter   $F_n$  has a density  $f_n(x) = \frac{dF_n(x)}{dx}$, $x \in [0,1]$, and a priori information is provided that there is positive constant  $C$ such that we have
\begin{equation*}
F_n \in \Xi(C) = \left\{F\,:\,\|f\|^2  < C,\, f(x) = \frac{dF(x)}{dx},\, F \in \Im \right\},
\end{equation*}
where $\|f\|^2 = \int_0^1 f^2(x) \,d\,x$.

Distribution function $F_n$ could be naturally replaced with distribution function $G_n$.

Main term of asymptotics for variance of chi-squared test statistics is significantly simplified if we suppose additionally 
\begin{equation*}
F_n \in \Xi_{1n} = \left\{F\,:\,\sup_{x \in [0,1]} |f(x)|  < c_n m_n^{1/2},\, f(x) = \frac{dF(x)}{dx},\, F \in \Im \right\},
\end{equation*}
where $c_n \to 0$ as $n \to \infty$.

For sequence $b_n > 0$,for $i=1,2$, define sets of alternatives $\Psi_i(b_n) = \{(F,G)\,:\, T_{in}(F-G) \ge b_n,\,\, (F,G) \in \Im\times\Im\}$.

We establish uniform consistency of test statistics  $T_{in}(\hat F_n-\hat G_{l_n})$, $i =1,2$, in problems of hypothesis testing 
\begin{equation*}
\mathbb{H}_0 \,\,:\,\, F_n(x) = G_{l_n}(x), \qquad x \in [0,1]
\end{equation*}
versus alternatives
\begin{equation*}
\mathbb{H}_n \,\,:\,\, (F_n,G_{l_n}) \in \Psi_i(b_n) \cap \Xi(C)
\end{equation*}
for sequences $b_n$, satisfying
\begin{equation}\label{i2}
0 < \liminf_{n \to \infty} m_n^{-1/2} b_n \le  \limsup_{n \to \infty} m_n^{-1/2} b_n <\infty.
\end{equation}
Suppose that  for all $j$, $1 \le j \le m_n$, we have
\begin{equation}\label{i3}
0 < c <\, m_n\, p_{nj}\,<C_1< \infty
\end{equation}
for some positive constants $c$ and $C_1$.

Suppose also that $m_n = o(n)$ as $ n \to \infty$.

Proof of Theorems is based on methods proposed in  \cite{er97} for the study of chi-squared tests for goodness-of-fit testing.

On set  $\Im \times \Im$ we define functional 
\begin{equation*}
T_n(F-G) = n \sum_{j=1}^m \left(\int_0^1 \phi_{nj}(x) d(F(x) - G(x))\right)^2 p_{nj}^{-1},
\end{equation*}
where $\phi_{nj}(x) = \mathbf{1}_{\{x \in I_{nj}\}} -  p_{nj}$, $x \in[0,1]$, $1 \le j \le m_n$ and $\mathbf{1}_{\{A\}}$-- denotes indicator of an event $A$. After that we explore test statistics as test statistics generated this functional.

This approach allows to prove easily the results similarly numerous results  \cite{er20,from12,ing02,tong}, established  for nonparametric hypothesis testing on a density based on expansions of series of orthogonal functions. However, in this case, functions $\phi_{nj}$ are not orthogonal.

In this notation test statistics $T_{2n}(\hat F_n-\hat G_{l_n})$ have the following form
\begin{equation*}
T_{2n}(\hat F_n-\hat G_{l_n}) = n \sum_{j=1}^n g_{nj}\left(\int_0^1 \phi_{nj}(x) d(\hat F_n(x) - \hat G_{l_n}(x))\right)^2 p_{nj}^{-1},
\end{equation*}
Note that, if hypothesis holds,  $\mathbf{E}[ T_{2n}(\hat F_n-\hat G_{l_n})]$ depends on unknown distribution function $F_n=G_{l_n}$, and, in the case of alternative,  $\mathbf{E}[ T_{2n}(\hat F_n-\hat G_{l_n})]$ depends on both unknown distribution functions $F_n$ and $G_{l_n}$. This is caused terms
\begin{equation*}
W_n = n \sum_{j=1}^n g_{nj}\left(\int_0^1 \phi_{nj}^2(x) d \hat F_n(x) +  \int_0^1 \phi_{nj}^2(x) d\,\hat G_{l_n}(x)\right) p_{nj}^{-1}
\end{equation*}
included in test statistics.

 To delete this dependence we subtract this term from test statistics in one of setups. Note  that we do not have such an influence of $W_n$ on test statistics  in the case of test statistics  $T_{1n}(\hat F_n-\hat G_{l_n})$.

Without loss of generality, we can suppose that distribution functions $F_n$ and $G_{l_n}$ have densities
\begin{equation*}
f_n(x) = 1 + \sum_{j=1}^m \theta_{nj}\,\phi_{nj}(x), \quad x \in [0,1]
\end{equation*}
and
\begin{equation*}
g_{l_n}(x) = 1 + \sum_{j=1}^m \tau_{nj}\,\phi_{nj}(x), \quad x \in [0,1]
\end{equation*}
respectively and
\begin{equation*}
\sum_{j=1}^m \theta_{nj} p_{nj} = 0, \qquad \sum_{j=1}^m \tau_{nj} p_{nj} = 0.
\end{equation*}
Denote $\eta_{nj} =  \theta_{nj} - \tau_{nj}$.

\subsection{Test statistics $T_{1n}$}
Denote $M_{1n}(\eta) = n m \sum_{j=1}^m p^2_{nj} \eta_{nj}^2$ and denote
\begin{equation*}\sigma_{1n}^2 = 2 m^2 \sum_{j=1}^m p_{nj}^2(1 + \theta_{nj} + a_n + a_n \tau_{nj})^2.
\end{equation*}
\begin{lemma}\label{l2a} We have
\begin{equation}\label{le2a1}
\mathbf{E} [T_{1n}(\hat F_n - \hat G_{l_n})] - (m-1) (1+ a_n)=   n m \sum_{j=1}^m p^2_{nj} \eta_{nj}^2(1 + o(1)),
\end{equation}
and
\begin{equation}\label{le2a3}
\begin{split}&
\mathbf{Var} [T_{1n}(\hat F_n - \hat G_{l_n})]
= \sigma_{1n}^2(1+ o(1))\\& + n m^2 \sum_{j=1}^m p_{nj}^3(1 + \theta_{nj} + a_n+ a_n \tau_{nj})\eta_{nj}^2(1+ o(1))\doteq \sigma_{11n}^2(1 + o(1))
\end{split}
\end{equation}
as $n \to \infty$.
\end{lemma}
Note that second addendum in right-hand side of (\ref{le2a3}) equal zero, if hypothesis holds. Thus we have interesting situation. The sets of alternatives is so reach that asymptotic variance for alternatives approaching to hypothesis is significantly different from asymptotic variance for hypothesis. 

By (\ref{lp36}), $\sigma_{11n}^2 - \sigma_{1n}^2 > 0$.
If $F_n \in \Xi_{2n}$, then $\sigma_{11n}^2 - \sigma_{1n}^2 = o(\sigma_{1n}^2)$ as $n \to \infty$.

Note that we can substitute into   (\ref{le2a3}) estimators
\begin{equation*}
\hat \theta_{nj} =\frac{\hat r_{nj}}{p_{nj}} - 1, \quad\mbox{and}\quad \hat\tau_{nj}= \frac{\hat s_{nj}}{p_{nj}} - 1
\end{equation*}
of parameters $\theta_{nj}$ and $\tau_{nj}$. After that, as we show, we get consistent estimator
\begin{equation*}
\hat \sigma^2_{1n} = 2m^2  \sum_{j=1}^m(\hat r_{nj} + a_n \hat s_{nj})^2
\end{equation*}
of variance $\sigma^2_{1n}$.

Other methods of estimation of variance are considered in \cite{and,from12,from13}.

Define tests
\begin{equation*}
K_{1n} = \mathbf{1}_{\{\hat\sigma_{1n}^{-1}(T_{1n}(\hat F_n-\hat G_{l_n}) - m(1 +a_n)) > x_\alpha\}},
\end{equation*}
where $x_{\alpha}$ is defined by equation $1 - \alpha = \Phi(x_\alpha)$, $0 < \alpha < 1$.
\begin{theorem} \label{th1} Assume (\ref{i2}), (\ref{i3})and let $m_n = o(n)$ as $ n \to \infty$. Then sequence of sets of alternatives $\Psi_{1n}(b_n) \cap \Xi(C)$ is uniformly consistent for sequence of tests  $K_{1n}$,generated tests statistics $ T_{1n}(\hat F_n-\hat G_{l_n})$.

We have  $\alpha(K_{1n}) = \alpha(1 +o(1))$ and
\begin{equation}\label{et11}
\beta(K_{1n},\Psi_{1n}(b_n)) = \Phi (\sigma_{11n}^{-1}(\sigma_{1n}x_\alpha - M_{1n}(\eta))) + o(1).
\end{equation}
as  $n \to \infty$.
\end{theorem}
\subsection{Test statistics $T_{2n}$ and $T_{3n}$} Denote $ M_{2n}(\eta) =  n \sum_{j=1}^m g_{nj}p_{nj} \eta_{nj}^2$ and denote
\begin{equation*}
\sigma_{2n}^2 = 2  \sum_{j=1}^m g^2_{nj}(1 + \theta_{nj} + a_n + a_n \tau_{nj})^2.
\end{equation*}
We show that 
\begin{equation*}
\hat \sigma^2_{2n} = 2  \sum_{j=1}^m g_{nj}^2 p_{nj}^{-2}(\hat r_{nj} + a_n \hat s_{nj})^2.
\end{equation*}
is consistent estimator of $\sigma_{2n}^2$.

Tests for test statistics   $ T_{2n}(\hat F_n-\hat G_{l_n})$ are based on the following asymptotics. 
\begin{lemma}\label{l2} We have
\begin{equation}\label{le21}
\mathbf{E} [T_{2n}(\hat F_n - \hat G_{l_n})] =  M_{2n}(\eta) (1 + o(1)) + \mathbf{E}[W_n] , \end{equation}
\begin{equation}\label{le21a}
\begin{split} &
\mathbf{E}[W_n] = \sum_{j=1}^m g_{nj} ((1 - p_{nj} + \theta_{nj}(1 - p_{nj}) - p_{nj} \theta_{nj}^2)\\& + a_n(1 - p_{nj} + \tau_{nj}(1 - p_{nj}) - p_{nj} \tau_{nj}^2))  \doteq e_n,
\end{split}
\end{equation}
\begin{equation}\label{le22}
\begin{split} &
\mathbf{Var} [T_{2n}(\hat F_n - \hat G_n)] = \sigma_{2n}^2(1+ o(1))\\& + n  \sum_{j=1}^m g_{nj}^2 p_{nj}(1 + \theta_{nj} + a_n+ a_n \tau_{nj})\eta_{nj}^2(1+ o(1))\doteq \sigma_{21n}^2(1 + o(1)).
\end{split}
\end{equation}
as $n \to \infty$.
\end{lemma}
As we show, if $m_n = o(n^{2/3})$, then
\begin{equation}\label{le24}
e_n = \sum_{j=1}^m g_{nj} (1 + a_n + \theta_{nj} + \tau_{nj}) + O(1).
\end{equation}
Note that we can substitute  estimators $\hat\theta_{nj}$ and $\hat\tau_{nj}$ of parameters $\theta_{nj}$ and  $\tau_{nj}$ into (\ref{le21a})  and to obtain consistent estimator  $\hat e_n$  for  $e_n$.

Define tests
\begin{equation*}
K_{2n} = \mathbf{1}_{\{\hat\sigma_{2n}^{-1}(T_{2n}(\hat F_n-\hat G_{l_n}) - \hat e_n) > x_\alpha\}},
\end{equation*}
where $x_{\alpha}$ is defined by equation $1 - \alpha = \Phi(x_\alpha)$, $0 < \alpha < 1$.
\begin{theorem} \label{th2} Assume (\ref{i2}),  (\ref{i3}) and let $m_n = o(n^{2/3})$ as $ n \to \infty$. Then sequence of sets of alternatives $\Psi_{2n}(b_n) \cap \Xi(C)$is uniformly consistent for sequence of tests $K_{2n}$.

Let $m_n = o(n)$ as $n \to \infty$ and let there be constant $C$ such that $\|g_n\| < C$, $g_n(x) = \frac{dG_n(x)}{dx}$, $x \in [0,1]$. Then sequence of sets of alternatives  $\Psi_{2n}(b_n) \cap \Xi(C)$ is uniformly consistent.

We have $\alpha(K_{2n}) = \alpha(1 +o(1))$ and
\begin{equation}\label{etua}
\beta(K_{2n},F_n,G_n) = \Phi (\sigma_{21n}^{-1}(\sigma_{2n}x_\alpha - M_{1n}(\eta))) + o(1)
\end{equation}
as $n \to \infty$.
\end{theorem}
In \cite{and,er97,from12,from13,tong} authors delete a version of addendum $W_n$ from version of test statistics  $T_{2n}$ for similar setups of nonparametric hypothesis testing and obtain the results for such modified test statistics.

Define test statistics
\begin{equation*}
T_{3n}(\hat F_n-\hat G_{l_n}) = T_{2n}(\hat F_n-\hat G_{l_n}) - W_n.
\end{equation*}
Define corresponding test of hypothesis testing
\begin{equation*}
K_{3n} = \mathbf{1}_{\{\hat\sigma_{2n}^{-1}T_{3n}(\hat F_n-\hat G_{l_n}) > x_\alpha\}},
\end{equation*}
where $x_{\alpha}$ is defined by equation $1 - \alpha = \Phi(x_\alpha)$, $0 < \alpha < 1$.
\begin{theorem} \label{th3} Assume (\ref{i2}),  (\ref{i3}) and let $m_n = o(n)$ as $ n \to \infty$. Then sequence of sets of alternatives $\Psi_{2n}(b_n) \cap \Xi(C)$ is uniformly consistent for sequence of tests $K_{3n}$.

 We have $\alpha(K_{3n}) = \alpha(1 +o(1))$ and
\begin{equation}\label{etub}
\beta(K_{3n},F_n,G_n) = \Phi (\sigma_{21n}^{-1}(\sigma_{2n}x_\alpha - M_{2n}(\eta))) + o(1).
\end{equation}
as $n \to \infty$.
\end{theorem}
\subsection{Hypothesis testing on homogeneity in terms of densities} 
Asymptotics of type II error probabilities in (\ref{et11}), (\ref{etua})and (\ref{etub}) for chi-squared tests or testing hypothesis of homogeneity are completely similar to asymptotics \cite{er97,er20} for goodness-of -fit testing of type II error probabilities of chi-squared tests for goodness-of -fit testing  (\ref{uuu}). By this reason, we can transfer necessary and sufficient conditions \cite{er20} of uniform consistency for chi-squared tests in the problem of goodness-of-fit testing to the case of hypothesis homogeneity. In \cite{er20} problem has been explored for alternatives defined in terms of densities.

Suppose distribution functions $F_n$ and $G_{l_n}$ have densities $f_n$, $g_{l_n}$ respectively and $F_n \in \Xi(C)$, $G_{l_n} \in \Xi(C)$. Denote $h_n = f_n - g_{l_n}$.

We explore problem of testing hypothesis
\begin{equation*}
\mathbb{H}_0\,:\, h_n(x) = 0, \qquad x \in [0,1],
\end{equation*}
versus alternatives
\begin{equation*}
\mathbb{H}_n\,:\, h_n \in \Omega_n \subset \Gamma,
\end{equation*}
where $\Gamma = \{h\,:\, h = \frac{ d(F- G)(x)}{d\,x}, \|h\| < \infty, F \in \Xi(C)\,\}$.

For this setup all statements of Theorem 6.1 in \cite{er20} hold if we replace densities  $1+f_n$ with functions $h_n$. All requirement in condition $B$ that functions  $1+f_n$ and functions specially defined by function $1+f_n$ were densities are replaced with the requirement that functions  $h_n$ and functions similarly specially defined by $h_n$ were differences of two densities.  In particular this holds if densities of distribution functions  $F_n$ and  $G_n$ satisfy $B$.
 
This version of Theorem 6.1 in \cite{er20} holds only for sequence of simple alternatives $h_n$, $\|h_n\| \asymp n^{-r}$, $\frac{1}{4} < r < \frac{1}{2}$, $m_n \asymp n ^{2 -4r}$. In this setup similarly to  \cite{er20}, we suppose that cells of chi-squared tests have the same length. 
\section{Proof of Theorems}
\subsection{Estimate of $\mathbf{E}\, [T_n]$} Reasoning will be provided for test statistics  $T_n$. Alternatives satisfy inequality
\begin{equation*}
T_n(F_n-G_n) = n \sum_{j=1}^m p_{nj} \eta_{nj}^2 \ge b_n.
\end{equation*}
By $f_n \in \Xi(C)$, we have 
\begin{equation}\label{o1}
\sum_{j=1}^m p_{nj}\theta_{nj}^2 \le \|f_n -1\|^2 < C.
\end{equation}
\begin{lemma}\label{l1} For $1 \le j \le m$ we have
\begin{equation}\label{le11}
\mathbf{E}_\theta [\phi_{nj}(X_1)] = \theta_{nj} p_{nj},
\end{equation}
\begin{equation}\label{le12}
\mathbf{E}_\theta [\phi_{nj}^2(X_1)] = p_{nj} (1 - p_{nj} + \theta_{nj}(1 - 2p_{nj})),
\end{equation}
\begin{equation}\label{le14}
\begin{split}&
\mathbf{E} [\bar \phi_{nj_1}^4(X_1)]= p_{nj}(1 + \theta_{nj})(1 - 4 p_{nj}(1 + \theta_{nj})\\& + 6 p^2_{nj}(1 + \theta_{nj})^2 - 3 p_{nj}^3(1 + \theta_{nj})^3)
\end{split}
\end{equation}
and, for $1 \le j_1 < j_2 \le m$, we have
\begin{equation}\label{le13}
\mathbf{E}_\theta [\phi_{nj_1}(X_1)\,\phi_{nj_2}(X_1)] = - p_{nj_1} p_{nj_2} (1  + \theta_{nj}(1 -2 p_{nj}) + \theta_{nj_2}(1 - 2p_{nj_2})),
\end{equation}
\begin{equation}\label{le15}
\begin{split}&
\mathbf{E}_\theta [\bar \phi_{nj_1}^2(X_1) \,\bar\phi_{nj_2}^2(X_1) ] \,=\, p_{nj_1} p_{nj_2} (1 + \theta_{nj_1})(1 + \theta_{nj_2})\\& \times (p_{nj_1}(1 + \theta_{nj_1}) + p_{nj_2}(1 + \theta_{nj_2}) - 3p_{nj_1}p_{nj_2}(1 + \theta_{nj_1})(1 + \theta_{nj_2})).
\end{split}
\end{equation}
\end{lemma}
Equalities (\ref{le11})--(\ref{le15}) are obtained  by straightforward calculations and proof is omitted. 

\begin{proof}[Proof of Lemma \ref{l2}] We begin with proof of (\ref{le21}). For $x,y \in [0,1]$, denote
\begin{equation*}
\bar\phi_{nj}(x) = \phi_{nj}(x) - \mathbf{E}_\theta \phi_{nj}(X_1) = \phi_{nj}(x) - \theta_{nj} p_{nj}
\end{equation*}
and
\begin{equation*}
\tilde\phi_{nj}(y) = \phi_{nj}(y) - \mathbf{E}_\tau \phi_{nj}(Y_1) = \phi_{nj}(y) - \tau_{nj} p_{nj}.
\end{equation*}
Then
\begin{equation}\label{lp21}
T_n(\hat F_n - \hat G_n) = I_{1n} + I_{2n} + I_{3n} + W_n,
\end{equation}
with
\begin{equation*}
I_{1n}  = 2\, I_{11n} + 2\,  I_{12n} + 2 \,I_{13n},
\end{equation*}
where
\begin{equation*}
I_{11n} = \sum_{1\le i_1 < i_2 \le n} U_{1n}(X_{i_1},X_{i_2}), \qquad I_{12n}= \sum_{1\le i_1 < i_2 \le l_n} U_{2n}(Y_{i_1},Y_{i_2})
\end{equation*}
and
\begin{equation*}
I_{13n}= \sum_{i_1 =1}^n \sum_{i_2=1}^{l_n} U_{3n}(X_{i_1},Y_{i_2}),
\end{equation*}
where
\begin{equation*}
\begin{split}&
U_{1n}(X_{i_1},X_{i_2}) = \sum_{j=1}^m \frac{\bar\phi_{nj}(X_{i_1}) \bar\phi_{nj} (X_{i_2})}{np_{nj}}, \\& U_{2n}(Y_{i_1},Y_{i_2}) = \sum_{j=1}^m \frac{\tilde\phi_{nj}(Y_{i_1}) \tilde\phi_{nj} (Y_{i_2})}{np_{nj}},
\end{split}
\end{equation*}
and
\begin{equation*}
U_{3n}(X_{i_1},Y_{i_2}) = \sum_{j=1}^m \frac{\bar\phi_{nj}(X_{i_1}) \tilde\phi_{nj} (Y_{i_2})}{np_{nj}}.
\end{equation*}
We have
\begin{equation}\label{lp22}
I_{2n} = \sum_{j=1}^m \left(\frac{1}{n} \sum_{i=1}^n \bar\phi_{nj}(X_i) - \frac{1}{l_n} \sum_{i=1}^{l_n} \tilde\phi_{nj}(Y_i)\right) \eta_{nj},
\end{equation}
\begin{equation}\label{lp23}
I_{3n} =  M_n(\eta) = n \sum_{j=1}^n p_{nj} \eta_{nj}^2 = T_n(F_n - G_{l_n}).
\end{equation}
\begin{equation}\label{lp24}
W_n = n^{-1}\sum_{j=1}^m  \sum_{i=1}^n \bar\phi_{nj}^2(X_i)p_{nj}^{-1} + n l_n^{-2}\sum_{j=1}^m  \sum_{i=1}^{l_n} \tilde\phi_{nj}^2(Y_i)p_{nj}^{-1}
\end{equation}
We have
\begin{equation}\label{lp24a}
\mathbf{E} I_{1n} = 0, \qquad \mathbf{E} I_{2n} = 0,
\end{equation}
\begin{equation}\label{lp25}
\begin{split}&
\mathbf{E} [W_{n}] = \sum_{j=1}^m (1 - p_{nj}  + \theta_{nj}(1 - 2p_{nj}) - p_{nj}\theta_{nj}^2)\\&
+ n l_n^{-1} \sum_{j=1}^m (1 - p_{nj}  + \tau_{nj}(1 - 2p_{nj}) - p_{nj}\tau_{nj}^2)\\&=
(1 + a_n)\sum_{j=1}^m (1 - p_{nj}  + \theta_{nj}(1 - 2p_{nj}) - p_{nj}\theta_{nj}^2)\\& + O(n^{-1/2} m M_{1n}^{1/2}(\eta))(1 + n^{-1}M_n(\eta))),
\end{split}
\end{equation}
because
\begin{equation}\label{lp26}
\begin{split}&
\sum_{j=1}^m |\theta_{nj} - \tau_{nj}| \le \max_{1 \le j \le m} p_{nj}^{-1} \sum_{j=1}^m p_{nj} |\eta_{nj}| \\& \le C m \left(\sum_{j=1}^m p_{nj} \eta_{nj}^2\right)^{1/2}   \left(\sum_{j=1}^m p_{nj}\right)^{1/2} \le C n^{-1/2} m M_{1n}^{1/2}(\eta)
\end{split}
\end{equation}
and
\begin{equation}\label{lp27}
\begin{split}&
\sum_{j=1}^m |\theta_{nj}^2 - \tau_{nj}^2| \le \max_{1 \le j \le m} p_{nj}^{-1} \sum_{j=1}^k p_{nj} |\eta_{nj}| \, (|\theta_{nj}| + |\tau_{nj}|)\\&  \le
C m^{-1} \left(\sum_{j=1}^m p_{nj} \eta_{nj}^2\right)^{1/2}   \left(\sum_{j=1}^m p_{nj}(\theta_{nj}^2 + \tau_{nj}^2)\right)^{1/2} \\&
\le  C n^{-1/2} m M_{1n}^{1/2}(\eta)  (N_n(\theta) + N_n(\tau))^{1/2} \le C n^{-1/2} m M_{1n}^{1/2}(\eta),
\end{split}
\end{equation}
because
\begin{equation*}
|N_n^{1/2}(\tau) - N_n^{1/2}(\theta)| \le n^{-1/2} M_n^{1/2}(\eta).
\end{equation*}
Note that reminder in right-hand side of (\ref{lp25}) is $o(m_n)$ as $n \to \infty$, if $m_n = o(n^{2/3})$.
\end{proof}
\subsection{Analysis of $\mathbf{Var}[T_n]$} We have
\begin{equation}\label{lp28}
\mathbf{Var}[I_{11n}] = V_{11n} + V_{12n},
\end{equation}
where
\begin{equation}\label{lp29}
\begin{split}&
V_{11n} = 2 \sum_{j=1}^m p_{nj}^{-2} (\mathbf{Var}[\phi_j(X_1)])^2\\&=
2 \sum_{j=1}^m (1 - p_{nj} + \theta_{nj}(1 -2p_{nj}) - p_{nj}^2)^2\\&=
2 \sum_{j=1}^m (1  + \theta_{nj})^2(1 +o(1))
\end{split}
\end{equation}
and
\begin{equation}\label{lp30}
\begin{split}&
V_{12n} =  2 \sum_{1 \le j_1 < j_2 \le m} p_{nj_1}^{-1} p_{nj_2}^{-1} (\mathbf{Cov}[\phi_{j_1}(X_1), \phi_{j_2}(X_1)])^2\\&=
4  \sum_{1 \le j_1 < j_2 \le m } p_{nj_1} p_{nj_2} (1 + \theta_{nj_1})^2(1 + \theta_{nj_2})^2(1 + o(1))\\& \le (C + N_n^2(\theta_n))(1 +o(1)).
\end{split}
\end{equation}
Therefore
\begin{equation}\label{lp30a}
\mathbf{Var}[I_{11n}] = 2 \sum_{j=1}^m (1 + \theta_{nj})^2(1 + o(1)).
\end{equation}
We have
\begin{equation}\label{lp31}
\begin{split}&
\mathbf{Var}[I_{12n}] =  4 a_n \sum_{j=1}^m p_{nj}^{-2} \mathbf{Var}[\phi_j(X_1)]  \mathbf{Var}[\phi_j(Y_1)]\\& =  4  a_n \sum_{j=1}^m (1 + \theta_{nj})(1 + \tau_{nj}) (1 + o(1)).
\end{split}
\end{equation}
Arguing similarly to $(\ref{lp30a})$, we get
\begin{equation}\label{lp30b}
\mathbf{Var}[I_{13n}] = 2 a_n^2\sum_{j=1}^m (1 + \tau_{nj})^2(1 + o(1)).
\end{equation}
We have
\begin{equation}\label{lp32}
\mathbf{Cov}[I_{11n},I_{12n}]= 0, \quad \mathbf{Cov}[I_{11n},I_{13n}]= 0, \quad \mathbf{Cov}[I_{12n},I_{13n}]= 0.
\end{equation}
Thus, by (\ref{lp30a})-- (\ref{lp32}), we get
\begin{equation}\label{lp33}
\mathbf{Var}[I_{1n}] =  2 \sum_{j=1}m (1 + a_n +\theta_{nj} + a_n \tau_{nj})^2(1 + o(1)).
\end{equation}
We have
\begin{equation}\label{lp34}
\mathbf{Var}[I_{2n}] = J_{21n} + J_{22n} + J_{23n} + J_{24n},
\end{equation}
with
\begin{equation}\label{lp35}
\begin{split}&
J_{21n} = 2 n^{-1}(n-1)^2  \sum_{1 \le j_1 < j_2 \le m} \mathbf{Cov}[\phi_{j_1}(X_1), \phi_{j_2}(X_1)] \eta_{nj_1}\eta_{nj_2}\\& =
 2 n^{-1}(n-1)^2  \sum_{1 \le j_1 < j_2 \le m} p_{nj_1} p_{nj_2} (1 + \theta_{nj_1})(1 + \theta_{nj_2})\eta_{nj_1}\eta_{nj_2}(1 + o(1))\\& \le C \left(\sum_{j=1}^m p_{nj}(1 + \theta_{nj_1})^2\right) \left(n \sum_{j=1}^m p_{nj} \eta_{nj}^2\right) \le C\,M_{1n}(\eta)(1 + N_n(\theta)),
\end{split}
\end{equation}
and
\begin{equation}\label{lp36}
\begin{split}&
J_{22n} = n^{-1} (n-1)^2 \sum_{j=1}^m  \mathbf{Var}[\phi_{nj}(X_1)] \,\eta_{nj}^2 \\& = n^{-1} (n-1)^2 \sum_{j=1}^m p_{nj} (1 - p_{nj} + \theta_{nj}(1 - 2p_{nj}) - p_{nj}\theta_{nj}^2) \,\eta_{nj}^2\\&=
 n\sum_{j=1}^m p_{nj} (1  + \theta_{nj}) \,\eta_{nj}^2(1 +o(1))=
 O(m^{1/2} M_{1n}(\eta)),
\end{split}
\end{equation}
because
\begin{equation}\label{lp41}
\max_{1 \le j \le m} |\theta_{nj}|^2 < C\, m\, N_{n}(\theta) < C\, m.
\end{equation}
Addendums $J_{23n}$ and $J_{24n}$ are estimated similarly to  $J_{21n}$ and $J_{22n}$ respectively. We omit this reasoning. 

We have
\begin{equation}\label{lp37}
\mathbf{Var}[W_n] = A_{1n} + A_{2n} + A_{3n} + A_{4n},
\end{equation}
where
\begin{equation}\label{lp38}
A_{1n} = n^{-1}  \sum_{1 \le j_1 < j_2 \le m} \mathbf{E}[\bar \phi_{nj_1}^2(X_1) \,\bar\phi_{nj_2}^2(X_1) ] p_{nj_1}^{-1} p_{nj_2}^{-1}
\end{equation}
and
\begin{equation}\label{lp39}
A_{2n} = n^{-1} \sum_{j=1}^m \mathbf{E} [\bar \phi_{nj_1}^4(X_1)]p_{nj}^{-2}.
\end{equation}
Addendums $A_{3n}$ and $A_{4n}$ are estimated similarly to  $A_{1n}$ and $A_{2n}$ respectively. We omit this reasoning.

Using (\ref{le14}) and (\ref{lp41}), we get
\begin{equation}\label{lp42}
\begin{split}&
A_{1n} \le n^{-1} \sum_{1 \le j_1 < j_2 \le m}[p_{nj_1} (1 + \theta_{nj_1})^2(1 + \theta_{nj_2})+ p_{nj_2}(1 + \theta_{nj_1}) (1 + \theta_{nj_2})^2] \\& \le Cn^{-1} \sum_{j=1}^m p_{nj}( 1 + |\theta_{nj}|)^2\left(m + \sum_{j=1}^m p_{nj} |\theta_{nj}|\right)\\& \le
Cn^{-1}( C + N_n(\theta))(m + m^{1/2}N^{1/2}(\theta)) \\&\le
Cn^{-1}m + Cn^{-1}m N_n(\theta) + C n^{-1} m^{1/2} N^{3/2}(\theta).
\end{split}
\end{equation}
Using (\ref{le15}) and (\ref{lp41}), we get
\begin{equation}\label{lp43}
\begin{split}&
A_{2n} =  n^{-1} \sum_{j=1}^m p_{nj}^{-1}(1 + \theta_{nj})[1 - 4 p_{nj}(1 + \theta_{nj})\\& + 6 p^2_{nj}(1 + \theta_{nj})^2 - 3 p_{nj}^3(1 + \theta_{nj})^3]
\end{split}
\end{equation}
We estimate only two addendums in  $A_{2n}$. Other two addendums are estimated similarly and have the smaller order.

We have
\begin{equation}\label{lp431}
\begin{split}&
 n^{-1} \sum_{j=1}^m p_{nj}^{-1}(1 + \theta_{nj}) \le C n^{-1}m^2\left(1 + \sum_{j=1}^m p_{nj}|\theta_{nj}|\right)\\&
\le C n^{-1}m^2\left(1 + \left(\sum_{j=1}^m p_{nj}\theta_{nj}^2\right)^{1/2}\right) \le C n^{-1}m^2(1 +N_n(\theta)) =o(m)
\end{split}
\end{equation}
and
\begin{equation}\label{lp432a}
\begin{split}&
n^{-1} \sum_{j=1}^m p_{nj}^2(1 + \theta_{nj})^4 \le Cn^{-1}m^{-1} + n^{-1} \sum_{j=1}^m p_{nj}^2 \theta_{nj}^4\\& \le Cn^{-1}(m^{-1}+ N_n^2(\theta))
\end{split}
\end{equation}
Therefore
\begin{equation}\label{lp432}
A_{2n} \le Cn^{-1}m^2(1 + N_n^{1/2}(\theta)) + n^{-1} N_n^2(\theta).
\end{equation}
\subsection{Consistency of estimators of bias and variance of test statistics $T_n$} Let us show consistency of estimators of   $\sum_{j=1}^m g_{nj}\theta_{nj}$ in  (\ref{le21a}) and (\ref{le22}).

We have
\begin{equation}\label{lp44}
\begin{split}&
 \mathbf{Var}\left[\sum_{j=1}^m g_{nj} \frac{1}{n} \sum_{i=1}^n \frac {\phi_{nj}(X_i)}{p_{nj}}\right] =
\frac{1}{n} \sum_{j=1}^m g_{nj}^2 \frac{\mathbf{Var}[\phi_{nj}(X_1)]}{p_{nj}^2}\\& + \frac{1}{n} \sum_{1 \le j_1 < j_2 \le m} g_{nj_1} g_{nj_2} \frac{\mathbf{Cov}[\phi_{nj_1}(X_1),\phi_{nj_2}(X_1)]}{p_{nj_1}p_{nj_2}}\\& = \frac{1}{n} \sum_{j=1}^m g_{nj}^2\frac{1 +  \theta_{nj}}{p_{nj}}(1 +o(1))\\& + \frac{1}{n} \sum_{1 \le j_1 < j_2 \le m} g_{nj_1} g_{nj_2} (1 + \theta_{nj_1} + \theta_{nj_2}) (1 + o(1)) = o(m),
\end{split}
\end{equation}
because
\begin{equation}\label{lp45}
n^{-1}\sum_{j=1}^m \frac{\theta_{nj}}{p_{nj} }\le Cn^{-1} m^2 \sum_{j=1}^m p_{nj} \theta_{nj} \le C n^{-1} m^2 N_n^{1/2}(\theta) = o(m)
\end{equation}
and
\begin{equation}\label{lp40}
\begin{split}&
n^{-1} m \sum_{j=1}^m g_{nj} \theta_{nj} \le C n^{-1} m \max_{1 \le j \le m} p_{nj}^{-1} \sum_{j=1}^m p_{nj} \theta_{nj}\\& \le C n^{-1} m^2 N_n^{1/2}(\theta) = o(m)
\end{split}
\end{equation}
We estimate only one addendums arising in the estimation of the variance. Other addendums are estimated similarly. 

We have
\begin{equation}\label{lp46}
n^{-4} \mathbf{Var}\left[\sum_{j=1}^m g^2_{nj} p_{nj}^{-2} \sum_{1  \le i_1 < i_2<n} \phi_{nj}(X_{i_1}) \phi_{nj}(X_{i_2})\right] \le B_{1n} + B_{2n},
\end{equation}
where
\begin{equation}\label{lp47}
\begin{split}&
B_{1n} = C n^{-2} \sum_{j=1}^m p_{nj}^{-4} (\mathbf{Var} [\phi_{nj}(X_1)])^2  \\& \le  C n^{-2} \sum_{j=1}^m p_{nj}^{-2}(1 +  \theta_{nj})^2(1 + o(1))\\& \le Cn^{-2} \max_{1 \le j \le m} p_{nj}^{-3} \sum_{j=1}^m p_{nj} (1 +\theta_{nj})^2 = O(n^{-2}m^3(1 + N_n(\theta)) = o(m)
\end{split}
\end{equation}
and
\begin{equation}\label{lp48}
\begin{split}&
B_{2n} = C n^{-2} \sum_{1 \le j_1 < j_2 \le m} \frac{(\mathbf{Cov}[\phi_{nj_1}(X_1),\phi_{nj_2}(X_1)])^2}{p_{nj_1}^2p_{nj_2}^2}\\& \le  C n^{-2} \sum_{1 \le j_1 < j_2 \le m}(1 +\theta_{nj_1} + \theta_{nj_2})^2 \le cn^{-2}m^2 \\& +
Cn^{-2}m \max_{1 \le j \le m} p_{nj}^{-1} \left(\left|\sum_{j=1}^m p_{nj} \theta_{nj}\right| + \sum_{j=1}^m p_{nj} \theta_{nj}^2\right)\\& \le Cn^{-2}m^2(1 + N_n(\theta)) = o(1).
\end{split}
\end{equation}
We provided estimates of variance in the case of sample $X_1,\ldots,X_n$. In the case of sample $Y_1,\ldots,Y_n$ we have different situation. In this case in final estimates  $N_n(\theta_n)$ isw replaces with $N_n(\tau_n)$. However in this case we can situation with  $N_n(\tau_n) \to \infty$ as $n \to \infty$.
However
\begin{equation}\label{lp49}
N_n^{1/2}(\tau_n) \le N_n^{1/2}(\theta_n) + n^{-1/2} M_n^{1/2}(\eta_n).
\end{equation}
Since $N_n^{1/2}(\theta_n) < C < \infty$, it suffices to show that, if, in final estimates, we replace $N_n^{1/2}(\theta_n)$ with  $n^{-1} M_n(\eta_n)$,  then these estimates will have smaller order than  $M_n^2(\eta_n)$.

Note that in  (\ref{lp25})--(\ref{lp48}) the largest orders in final estimates for distribution function   $G_{l_n}$ are  $M_{1n}(\eta_n) N_n(\tau_n)$ (version of(\ref{lp35})), $n^{-1}m^2 N_n^{1/2}(\tau_n)$ (version of (\ref{lp42})) and $n^{-1} N_n^2(\tau_n)$ (version of (\ref{lp432})).

It suffices to estimate only   $n^{-1}m^2 N_n^{1/2}(\tau_n)$. We have
\begin{equation}\label{lp50}
n^{-3/2}m^2 M_n^{1/2}(\eta_n) M_n^{-2}(\eta_n) = O(n^{-3/2}m_n^2 m_n^{-3/4}) = o(1),
\end{equation}
if $m_n^{-1/2}M_n(\eta_n) \to \infty$ as $n \to \infty$.

Thus
\begin{equation}\label{lp51}
M_n(\eta_n) \hat \sigma_n \to_P \infty,
\end{equation}
if $m_n^{-1/2}M_n(\eta_n) \to \infty$ as  $n \to \infty$.

Therefore type II error probabilities of tests $K_n$ tends to zero if
$N_n(\tau_n) \to \infty$ as $n \to \infty$.

\subsection{Asymptotic normality of test statistics $T_n$} It suffices to  prove asymptotic normality of statistics $I_{1n}$. For alternatives we can suppose  $(F_n,G_{l_n}) \in \Xi_n(C) \times \Xi_n(C)$ for some $C > 0$.  Otherwise, type II error probabilities tends to zero. Statistics $I_{1n}$ are not $U$--statistics. However we can implement the same martingale technique to the proof of asymptotic normality   \cite{br,h,er97,ing02} and to get similar result as in the case of goodness-of-fit tests \cite{er97,ing02}. Since in \cite{and} similar reasoning for testing of hypothesis of homogeneity are omitted for test statistics based on  $\mathbb{L}_2$ --norm of kernel estimator of density we outline this reasoning in this paper.

The reasoning will be provided if  $l_n \le n$. The case $l_n \ge n$ is similar.

Define martingale $W_{ni}$, $1 \le i \le n + l_n$, by induction. We put
\begin{equation*}
W_{n1}  = U_{1n}(X_1,X_1),\quad \mbox{and} \quad W_{n2} = U_{2n}(Y_1,Y_1) + U_{3n}(X_1,Y_1).
\end{equation*}
If $i$ is odd,we put $j = [i/2]$ and
\begin{equation*}
W_{ni}  =  \sum_{s=1}^{j}  U_{1n}(X_j,X_s) + \sum_{s=1}^{j-1}  U_{3n}(X_j,Y_s)
\end{equation*}
If  $i$is even, $i \le 2l_n$, we put $j =i/2$ and
\begin{equation*}
W_{ni}  =  \sum_{s=1}^{j}  U_{2n}(Y_j,Y_s) + \sum_{s=1}^{j-1}  U_{3n}(X_s,Y_j)
\end{equation*}
If $i \ge 2l_n$, we put $j = i -l_n$ and
\begin{equation*}
W_{ni}  =  \sum_{s=1}^{j}  U_{1n}(X_j,X_s) + \sum_{s=1}^{l_n}  U_{3n}(X_j,Y_s).
\end{equation*}
We can implement to this martingale the reasoning of  \cite{h} and obtain similar result. 

Denote
\begin{equation*}
\begin{split}&
V_{1n}(x,y) = \mathbf{E}[U_{1n}(x,X_1) U_{1n}(y,X_1)],\,\,\, V_{2n}(x,y) = \mathbf{E}[U_{1n}(x,Y_1) U_{1n}(y,Y_1)],\\&
V_{3n}(x,y) = \mathbf{E}[U_{3n}(X_1,x) U_{3n}(X_1,y)], \,\,\, V_{4n}(x,y) = \mathbf{E}[U_{3n}(x,Y_1) U_{3n}(y,Y_1)].
\end{split}
\end{equation*}
\begin{theorem} \label{tm1}
Statistics $I_{1n}$ is asymptotically normal with zero mean and variance  $\sigma_1^2$, if we have
\begin{equation} \label{mor}
\begin{split}&
\lim_{n \to \infty} m_n^{-1} [ \mathbf{E}[V_{1n}^2(X_1,X_2) + V_{2n}^2(Y_1,Y_2) +V_{3n}^2(X_1,X_2) + V_{4n}^2(Y_1,Y_2)]\\& + n^{-1} \mathbf{E}[U_{1n}^4(X_1,X_2) + U_{2n}^4(Y_1,Y_2) + U_{3n}^4(X_1,Y_1)]] =0.
\end{split}
\end{equation}
\end{theorem}
Proof of Theorem almost repeat the reasoning for the proof of asymptotic normality in \cite{h} and is omitted.

 Checking (\ref{mor})  does not differ practically from checking similar conditions in the case of goodness-of-fit testing  \cite{er97}. Moreover the most part of estimates for proof of  (\ref{mor}) and estimates in \cite{er97} is coincide. Thus we omit this reasoning.

\end{document}